\documentclass[%
 aip,
cp,  % Conference Proceedings
 amsmath,amssymb,%nobibnotes,
 reprint,%
]{revtex4-2}

\usepackage{graphicx}% Include figure files
\usepackage{caption}
\usepackage{subcaption}
\usepackage{dcolumn}% Align table columns on decimal point
\usepackage{bm}% bold math
\usepackage[utf8]{inputenc}
\usepackage[T1]{fontenc}
\usepackage{mathptmx}

\newtheorem{theorem}{\quad Theorem}[section]
\newtheorem{prop}[theorem]{\quad Proposition}

\newtheorem{definition}{\quad Definition}[section]
\newtheorem{remark}[theorem]{\quad Remark}

\begin{document}

\title{Focal Curves of Closed  Toroidal Curves}% Force line breaks with \\
\author{Dinkova C. L.} % Write as First name Surname
 \email[Corresponding author: Cvetelina Lachezarova Dinkova ]{c.dinkova@shu.bg}
\author{Encheva R. P.}%
 \email[Radostina Petrova Encheva ]{r.encheva@shu.bg}
\affiliation{
  Faculty of Mathematics and Informatics, Konstantin Preslavsky University, Universitetska Str. 115\\ 9700 Shumen, Bulgaria% Force line breaks with \\ if necessary
}

\author{Ali A. A.}
 \email[Aynur Abdulova Ali ]{a.ali@shu.bg}
%\affiliation{%
%Second institution and/or address% Force line breaks with \\ if necessary
%}%
%\affiliation{You would list an author's second affiliation (if applicable) here.}

\date{\today} % It is always \today, today, but any date may be explicitly specified
              % Not printed for conference proceedings

\begin{abstract}
Geometric constructions are widely used in computer graphics and engineering drawing. A right generalized
cylinder is a ruled surface whose base curve is a plane curve perpendicular to the rulings. The paper discusses relations between the base curve and the non-planar space curve on the right generalized cylinder. Based on these relations, a method for obtaining a new space curve from a given plane curve parameterized about an arbitrary parameter is presented. At first, we define a non-planar space curve on the right generalized cylinder whose base curve is the considered plane curve,
parameterized about an arbitrary parameter. Later on, we examine the focal curve of the obtained cylindrical curve which is also a non-planar curve. The Frenet-Seret system of the cylindrical curve and its focal curve are expressed in terms of the signed curvature of the abovementioned plane curve and its derivatives. Finally, we obtained the parametric representation of orthogonal projection of the focal curve onto the Euclidean plane via before mentioned plane curve and its derivatives.
That curve is called generalized focal curve of a plane curve. The proposed method is
demonstrated for several closed plane curves used in engineering practice. These curves include: epicycloid, hypocycloid and a curve that is orthogonal projection of toroidal helix onto the Euclidean plane.
\end{abstract}


\maketitle

\section{\label{sec:level1}Introduction}

A closed curve is a mathematical concept used in geometry and topology. It is a continuous and connected curve that begins and ends at the same point. In other words, a simple closed curve is a curve that returns to its starting point without intersecting itself, whereas a complex closed curve intersects itself at one or more points. Closed curves can have various shapes and forms, and they are fundamental in many areas of mathematics and science, including physics, engineering, and computer graphics.
One of the most basic examples of a closed curve is a circle. A circle is a curve that is defined by a set of points that are equidistant from a central point. The circle is a simple closed curve because it forms a complete loop, and it is continuous because it has no gaps or breaks. Other examples of closed curves include ellipses, squares, rectangles, and polygons. These curves can have more complex shapes, and they are often defined by mathematical formulas or equations.
A cardioid is another simple closed curve that is defined as the set of all points in a plane that are at a fixed distance from a given point, called the focus, and whose paths are traced by a point on a circle rolling around a fixed circle. The cardioid has many important applications in mathematics, physics, and engineering.
One important characteristic of a closed curve is its perimeter or circumference. The perimeter is the total length of the curve, and it is measured in units such as meters, centimeters, or inches. For a circle, the circumference is given by the formula $C = 2\pi r$, where r is the radius of the circle and $\pi$ is the mathematical constant pi. For other closed curves the perimeter can be calculated using various formulas and methods, depending on the shape and size of the curve.
Closed curves also have a topological property known as orientability. An orientable closed curve is a curve that has a consistent orientation, which means that the curve can be traced in a single direction without crossing over itself. For example, a circle is an orientable curve because it has a consistent clockwise or counterclockwise direction. However, some closed curves are non-orientable, which means that they cannot be consistently traced in a single direction without crossing over themselves. An example of a non-orientable closed curve is the Möbius strip, which is a twisted loop that has only one side and one edge.

Closed curves are used in many practical applications, such as in the design of buildings, bridges, and roads. For example, architects and engineers use closed curves to design arches, domes, and other curved structures that are aesthetically pleasing and structurally sound. Closed curves are also used in computer graphics and animation to create smooth and continuous shapes and movements. In addition, closed curves are used in physics and engineering to model the behavior of fluids, electromagnetic fields and other physical phenomena.

Geometric constructions are widely used in computer graphics and engineering drawing. A right generalized
cylinder is a ruled surface whose base curve is a plane curve perpendicular to the rulings. In this paper, relations between
the base curve and the non-planar regular curve on the right generalized cylinder are discussed. Based on these relations, a
method for obtaining a new space curve from a given regular plane curve parameterized about an arbitrary parameter is presented. Furthermore a new curve associated with given plane curve, parameterized about an arbitrary parameter is defined.

In this research we consider examples of two popular classes of plane curves known as epicycloid and hypocycloid. An \textbf{epicycloid} is a plane curve produced by tracing the path of a chosen point on the circumference of a circle called an epicycle, which rolls without slipping around a fixed circle. An \textbf{hypocycloid} is a plane curve traced by a point on the circumference of a circle rolling internally on the circumference of a fixed circle.
An epicycloid (an hypocycloid) and its evolute (a focal curve) are similar.
In addition we examine another popular space curve called \textbf{toroidal helix} and its projection onto Euclidean plane.

Initially, we proposed a method for obtaining a new space curve lying on a torus from a given plane curve. After that we apply this method for abovementioned closed plane curves. Then we calculate the focal curves of obtained space curve which is an intersection curve of right generalized cylinder and a torus. Finally, we received a generalized focal curve of closed plane which is totally different from its evolute (a focal curve).

\section{Preliminaries}

This section introduces some basic concepts of the classical differential geometry of curves and surfaces in two and three dimensional Euclidean space. More details can be found in "Modern Differential Geometry of Curves and Surfaces" (see \cite{Gray2006}) and "Differential Geometry of Curves and Surfaces" (see \cite{Do Carmo2016}).

\begin{definition}\cite[p.18]{Do Carmo2016}
  A parameterized differentiable curve is a differentiable
map $\alpha: I \rightarrow \mathbb{E}^3$ of an open interval $I\subseteq \mathbb{R}$ of the real line $\mathbb{R}$ into $\mathbb{E}^3$.
\end{definition}
The word differentiable in this definition means that $\alpha$ is a correspondence
which maps each $t \in I$ into a point $\alpha(t) = (x(t), y(t), z(t))\in \mathbb{E}^3$,
in such a way that the functions $x(t), y(t), z(t)$ are differentiable. The variable $t$ is called
the parameter of the curve.

\begin{definition}\cite[p.22]{Do Carmo2016}
 A parameterized differentiable curve $\alpha: I \rightarrow \mathbb{E}^3$
is said to be regular if $\dot{\alpha}(t)=\frac{d\alpha(t)}{dt}\neq 0$ for all $t \in I$.

Given $t_0 \in I$ , the arc length of a regular parameterized curve $\alpha: I \rightarrow \mathbb{E}^3$,
from the point $t_0$, is by definition $s(t) =\int_{t_0}^{t} \|\dot{\alpha}(t)\| dt$,
where $\|\dot{\alpha}(t)\|=\sqrt{\dot{x}^{2}(t)+\dot{y}^{2}(t)+\dot{z}^2(t)}$
is the length of the vector $\dot{\alpha}(t)$. Since $\dot{\alpha}(t)\neq 0$, the arc length $s(t)$ is a
differentiable function of $t$ and $ds/dt = \|\dot{\alpha}(t)\|$.

A map $\alpha: I \rightarrow \mathbb{E}^3$
is called a curve of class $C^k$
if each of the coordinate functions in the expression $\alpha(t) = (x(t), y(t), z(t))$ has continuous
derivatives up to order $k$. If $\alpha$ is merely continuous, we say that $\alpha$ is of
class $C^0$. A curve $\alpha$ is called simple if the map $\alpha$ is one-to-one.

\end{definition}
\subsection{\label{sec:level2}Curves in Euclidean 2-space}

The image of any parameterized curve in the Euclidean plane $\mathbb{E}^2\equiv O\vec{\mathbf{e}}_1\vec{\mathbf{e}}_2$ under an orientation-preserving affine map in $\mathbb{E}^2$ is also a parameterized curve in $\mathbb{E}^2$. Now we will discuss the differential-geometric invariants of regular plane curves with respect to the group of orientation-preserving rigid motions.

A closed plane curve is a regular parameterized curve $\alpha: [a,b] \rightarrow \mathbb{E}^2$
such that $\alpha$ and all its derivatives agree at $a$ and $b$, that is
$\dot{\alpha}(a)=\dot{\alpha}(b),\, \ddot{\alpha}(a)=\ddot{\alpha}(b),\, \dddot{\alpha}(a)=\dddot{\alpha}(b),\,\ldots$

The curve $\alpha$ is simple if it has no further self-intersections, that is, if $t_1, t_2 \in
[a, b], t_1 \neq t_2$, then $\alpha(t_1) \neq \alpha(t_2)$.

Consider a regular plane curve $\mathbf{\alpha}:I \rightarrow \mathbf{E}^2$ of class $C^3$ that is defined on the open interval $I\subseteq \mathbb{R}$ by
\begin{equation}\label{plane_curve}
  \alpha:\alpha(t)=(x(t),y(t),0).
\end{equation}

\begin{definition}\cite[p.3]{Gray2006}
  The \textbf{complex structure} is the linear map $J:\mathbb{E}^2\rightarrow \mathbb{E}^2$ given by $J(p_1,p_2)=(-p_2,p_1).$ Geometrically, $J$ is rotation by $\pi/2$ in a counterclockwise direction.
\end{definition}

 We use the symbols $\dot{\alpha}=\displaystyle\frac{d\alpha(t)}{dt}$, $\ddot{\alpha}=\displaystyle\frac{d\dot{\alpha}(t)}{dt}$ for differentiation about an arbitrary parameter $t$.
The scalar product of two vector functions $x(t)=(x_1(t),x_2(t))$ and $y(t)=(y_1(t),y_2(t))$ is given by $<x(t),y(t)>=x_1(t)y_1(t)+x_2(t)y_2(t)$.
The norm of vector function $x(t)$ is given by $\|x(t)\|=\sqrt{<x,x>}=\sqrt{x_{1}^{2}(t)+x_{2}^{2}(t)}$ for $t\in \mathbb{R}$.

\begin{definition}\label{def.eucl.curv.}\cite[p.11]{Gray2006}
  Let $\mathbf{\alpha}:I\rightarrow\mathbb{E}^2$ be a regular curve. The \textbf{signed curvature} $K$ of $\mathbf{\alpha}$ is given by the formula $K(t)=\displaystyle\frac{<\mathbf{\ddot{\alpha}}(t),J\mathbf{\dot{\alpha}}(t)>}{\|\mathbf{\dot{\alpha}}(t)\|^3}$. The function $R=\displaystyle\frac{1}{K}$ is called the \textbf{radius of curvature} of $\mathbf{\alpha}.$
\end{definition}

\begin{remark}
The signed curvature $K$ defined by the equation in Definition \ref{def.eucl.curv.} above is invariant under orientation-preserving rigid motions in $\mathbb{E}^2$.
\end{remark}

\subsection{\label{sec:level3}Curves in Euclidean 3-space}

The Frenet-Seret system of a regular space curve $\gamma$ consists of a vector and scalar invariants of a curve:
 %$\gamma(t)=\alpha(t)+f(t)\vec{\mathbf{e}}_3$ a cylindrical curve over a plane curve $\alpha$;
 \[T=\displaystyle\frac{\dot{\gamma}}{\|\dot{\gamma}\|},\,
 N=\displaystyle\frac{(\dot{\gamma}\times\ddot{\gamma})\times \dot{\gamma}}
 {\|(\dot{\gamma}\times\ddot{\gamma})\times \dot{\gamma}\|},\,
 B=\displaystyle\frac{\dot{\gamma}\times\ddot{\gamma}}
 {\|\dot{\gamma}\times\ddot{\gamma}\|},\,\varkappa=\displaystyle\frac{\|\dot{\gamma}\times\ddot{\gamma}\|}{\|\dot{\gamma}\|^{3}},\,
 \tau=\displaystyle\frac{\dot{\gamma}\ddot{\gamma}\dddot{\gamma}}{\|\dot{\gamma}\times\ddot{\gamma}\|^{2}},\]
where vector invariants $T, N, B$ called a tangent, a principal and a binormal unit vector fields of $\gamma$, scalar invariants $\varkappa$ and $\tau$ called curvature and torsion of $\gamma$.
% $c_{1}(t)=\displaystyle\frac{1}{\varkappa(t)}$, $c_{2}(t)=-\displaystyle\frac{\frac{d}{dt}\varkappa(t)}{\parallel \frac{d\gamma(t)}{dt}\parallel \varkappa(t)^2 \tau(t)}=\displaystyle\frac{\frac{dc_1(t)}{dt}}{\parallel\frac{d\gamma(t)}{dt}\parallel \tau(t)}$ focal curvatures of $\gamma$.
% The focal curve of an immersed smooth curve $\gamma$, in Euclidean space $\mathbb{R}^{3}$, consists of the centres of its osculating spheres
%    \begin{equation}\label{foc}
%C_{\gamma}(t)=\gamma+c_1 N+c_2 B.
%    \end{equation}

\begin{definition}\cite[p.241]{Gray2006}
  The \textbf{focal curve} of a regular $C^3$ space curve $\mathbf{\gamma}: I\rightarrow\mathbb{E}^3$ is the curve given by
  \begin{equation}\label{C_gamma}
   \mathbf{ C_\gamma}(t)=\mathbf{\gamma}(t)+c_1(t) \mathbf{N}(t)+c_2(t) \mathbf{B}(t),
  \end{equation}
  where $\mathbf{N}$ is a principal unit normal vector field of $\mathbf{\gamma}$, $\mathbf{B}$ is a binormal unit vector field of $\mathbf{\gamma}.$ The coefficients $c_1(t)$ and $c_2(t)$ are smooth functions called focal curvatures of $\mathbf{\gamma}$, and given by
  \begin{equation}\label{c1,c2}
    c_1(t)=\frac{1}{\varkappa(t)},\quad c_2(t)=-\frac{\frac{d}{dt}\varkappa(t)}{\parallel \frac{d\mathbf{\gamma}(t)}{dt}\parallel \varkappa(t)^2 \tau(t)}=\frac{\frac{dc_1(t)}{dt}}{\parallel\frac{d\mathbf{\gamma}(t)}{dt}\parallel \tau(t)},
  \end{equation}
  where $\varkappa(t)$ and $\tau(t)$ are the Euclidean curvatures of $\mathbf{\gamma}$.
\end{definition}

 In other words, the focal curve of an immersed smooth curve $\gamma$, in Euclidean space $\mathbb{E}^{3}$, consists of the centres of its osculating spheres.

\begin{remark}
The functions $c_1(t)$ and $c_2(t)$ are well defined because $\varkappa(t)$ and $\tau(t)$ are non-zero functions.
\end{remark}

\subsection{\label{sec:level4}Surfaces in Euclidean 3-space}
\begin{definition}\cite[p.438]{Gray2006}
  Let $S\subset\mathbb{E}^3$ be a surface. Then $S$ is a \textbf{generalized cylinder} over a curve $\mathbf{\alpha} : I\rightarrow\mathbb{E}^3$ if  $S$ can be parameterized as
  \begin{equation}\label{gen:cyl}
  S: S(u,v)=\mathbf{\alpha}(u)+v.\vec{\mathbf{q}},
  \end{equation}
  where $\vec{\mathbf{q}}\in\mathbb{E}^3$ is a fixed vector.
\end{definition}

We consider a case of right generalized cylinder over the plane curve $\alpha$, when its rulings are perpendicular to the generating plane curve.
That means, the fixed vector is a unit vector $\vec{\mathbf{e}}_3=(0,0,1)$ $\parallel Oz$ and parameter $v$ is replaced by function
$f(v)\in \mathbb{R},\, f(v)\in C^{3}$. Then parametrisation of a \textbf{right generalized cylinder} is
\begin{equation}\label{right_gen_cyl}
  S_{1}(u,v)=\alpha(u)+f(v)\vec{\mathbf{e}}_3.
\end{equation}

\section{\label{sec:level5}Previous results}

\begin{theorem}\cite{Georgiev2015}
Let $\alpha(t)=(x(t),y(t),0), \, t\in I\subset \mathbb{R}$ be regular plane curve of class $C^3$ with a
 nonzero signed curvature, and let $f(t)\in C^3$ be a real-valued function. Suppose that $\vec{\mathbf{e}}_3$ is
  the unit vector on $Oz$-axis and
  \[
  \gamma(t)=\alpha(t)+f(t)\vec{\mathbf{e}}_3, \quad t\in I
  \]
  is a parameterized space curve. Then, $\gamma(t)$
   is a regular curve whose curvature and torsion are given by
   \begin{equation}\label{GENC}
   \varkappa=\frac{\sqrt{\big\langle \ddot{\alpha},J\dot{\alpha}\big\rangle^2+
   \big\langle \ddot{f}\dot{\alpha}-\dot{f}\ddot{\alpha},\ddot{f}\dot{\alpha}
   -\dot{f}\ddot{\alpha}\big\rangle}}{\left(\sqrt{\big\langle \dot{\alpha},
   \dot{\alpha}\big\rangle +\dot{f}^2}\right)^3}
   \end{equation}
\begin{equation}\label{GENT}
\tau=\frac{\dddot{f}\big\langle \ddot{\alpha},J\dot{\alpha}\big\rangle+\ddot{f}\big\langle
-J\dot{\alpha},\dddot{\alpha}\big\rangle
+\dot{f}<\big\langle J\ddot{\alpha},\dddot{\alpha}\big\rangle}
 {\big\langle \ddot{\alpha},J\dot{\alpha}\big\rangle^2+
   \big\langle \ddot{f}\dot{\alpha}-\dot{f}\ddot{\alpha},\ddot{f}\dot{\alpha}
   -\dot{f}\ddot{\alpha}\big\rangle}
 \end{equation}
 \end{theorem}

\section{\label{sec:level6}Research results}

The next statement gives us the relations between the Frenet-Seret frame of $\mathbf{\gamma}$ and the signed curvature as well as the arc-length prime of $\mathbf{\alpha}$, parameterized about an arbitrary parameter $t$.
\begin{theorem}\label{frenet frame via base curve}
  Let $\mathbf{\alpha}=\mathbf{\alpha}(t),\, t\in I\subseteq \mathbb{R}$ be a regular $C^3$-plane curve in $\mathbb{E}^2$ and $K\neq0$ is the Euclidean signed curvature of $\mathbf{\alpha}$. Let $\mathbf{\gamma}(t)=\mathbf{\alpha}(t)+f(t).\vec{\mathbf{e}}_3,\,t\in I\subseteq\mathbb{R},\, f(t)\in C^2$ be the corresponding cylindrical curve over the right generalized cylinder with a base curve $\mathbf{\alpha}.$ If $T,N,B$ are vector invariants of $\mathbf{\gamma}$, then they can be expressed by the derivatives of $\alpha$, the scalar function $f$ and the unit vector $\vec{\mathbf{e}}_3$ with the following equations
\begin{equation}\label{T}
  T=\displaystyle\frac{\dot{s}+\dot{f}\vec{e}_{3}}{\sqrt{\dot{s}^{2}+\dot{f}^{2}}}
\end{equation}
\begin{equation}\label{N}
  N=\displaystyle\frac{(\dot{s}^{2}+\dot{f}^{2})\ddot{\alpha}
 -\frac{1}{2}\frac{d}{dt}(\dot{s}^{2}+\dot{f}^{2})\dot{\alpha}+
 \left(\ddot{f}\dot{s}^{2}-\dot{f}\frac{d}{dt}\left(\frac{\dot{s}^{2}}{2}\right)\right)\vec{e}_{3}}
 {\sqrt{\dot{s}^{2}+\dot{f}^{2}}\,\sqrt{\dot{s}^{6}K^{2}+\|\ddot{f}\dot{\alpha}-\dot{f}\ddot{\alpha}\|^{2}}}
\end{equation}
\begin{equation}\label{B}
  B=\displaystyle\frac{-J(\ddot{f}\dot{\alpha}-\dot{f}\ddot{\alpha})+\dot{s}^{6}K^{2}\vec{e}_{3}}
 {\sqrt{\dot{s}^{6}K^{2}+\|\ddot{f}\dot{\alpha}-\dot{f}\ddot{\alpha}\|^{2}}}
\end{equation}
where $\dot{s}$ is the arc-length prime of $\mathbf{\alpha}$ with respect to arbitrary parameter $t.$
\end{theorem}
{\it Proof:} We apply a differentiation about arbitrary parameter $t$. From $\mathbf{\gamma}(t)=\mathbf{\alpha}(t)+f(t).\vec{\mathbf{e}}_3$ it follows that
\[\mathbf{\dot{\gamma}}=\frac{d\mathbf{\gamma}}{dt}=\mathbf{\dot{\alpha}}+\dot{f}.\vec{\mathbf{e}}_3,\quad
 \mathbf{\ddot{\gamma}}=\frac{d^2\gamma}{dt^2}=\mathbf{\ddot{\alpha}}+\ddot{f}.\vec{\mathbf{e}}_3 \quad\textrm{and}\quad
 \mathbf{\dddot{\gamma}}=\frac{d^3\gamma}{dt^3}=\mathbf{\dddot{\alpha}}+\dddot{f}.\vec{\mathbf{e}}_3.\]
Since $\mathbf{\dot{\alpha}}$, $J(\mathbf{\dot{\alpha}})$, $\vec{\mathbf{e}}_3$ and $\mathbf{\ddot{\alpha}}$, $J(\mathbf{\ddot{\alpha}})$, $\vec{\mathbf{e}}_3$ forms a right-handed orthogonal frame of the curve $\mathbf{\alpha}$ and cross vector product $\mathbf{\dot{\alpha}}\times\mathbf{\ddot{\alpha}}=<\mathbf{\mathbf{\ddot{\alpha}}},J\mathbf{\dot{\alpha}}>.\vec{\mathbf{e}}_3=\dot{s}^3K.\vec{\mathbf{e}}_3$, then for cross vector product of the first and second derivatives of $\mathbf{\gamma}$ we have
 \[\mathbf{\dot{\gamma}}\times\mathbf{\ddot{\gamma}}=(\mathbf{\dot{\alpha}}+\dot{f}.\vec{\mathbf{e}}_3)\times(\mathbf{\ddot{\alpha}}+\ddot{f}.\vec{\mathbf{e}}_3)=
 \mathbf{\dot{\alpha}}\times\mathbf{\ddot{\alpha}}-J(\ddot{f}\mathbf{\dot{\alpha}}-\dot{f}\mathbf{\ddot{\alpha}})=
 <\mathbf{\mathbf{\ddot{\alpha}}},J\mathbf{\dot{\alpha}}>.\vec{\mathbf{e}}_3-J(\ddot{f}\mathbf{\dot{\alpha}}-\dot{f}\mathbf{\ddot{\alpha}})=
 \dot{s}^3K.\vec{\mathbf{e}}_3-J(\ddot{f}\mathbf{\dot{\alpha}}-\dot{f}\mathbf{\ddot{\alpha}}).\]
 Then the norms of tangent vector $\mathbf{\dot{\gamma}}$ and binormal vector $\mathbf{\dot{\gamma}}\times \mathbf{\ddot{\gamma}}$ of $\mathbf{\gamma}$ are given by
 \[\|\mathbf{\dot{\gamma}}\|=\sqrt{<\mathbf{\dot{\gamma}},\mathbf{\dot{\gamma}}>}=\sqrt{<\mathbf{\dot{\alpha}},\mathbf{\dot{\alpha}}>+\dot{f}^2}=\sqrt{\dot{s}^2+\dot{f}^2},\quad
 \|\mathbf{\dot{\gamma}}\times \mathbf{\ddot{\gamma}}|
% = \sqrt{< \dot{s}^3K.\vec{\mathbf{e}}_3-J(\ddot{f}(t)\mathbf{\dot{\alpha}}(t)-\dot{f}(t)\mathbf{\ddot{\alpha}}(t)), \dot{s}^3K.\vec{\mathbf{e}}_3-J(\ddot{f}(t)\mathbf{\dot{\alpha}}(t)-\dot{f}(t)\mathbf{\ddot{\alpha}}(t))>}
 =\sqrt{\dot{s}^{6}K^{2}+\|\ddot{f}\dot{\alpha}-\dot{f}\ddot{\alpha}\|^{2}}\,.\]
 Finally, from $T=\displaystyle\frac{\mathbf{\dot{\gamma}}}{\|\mathbf{\dot{\gamma}}\|}$, $B=\displaystyle\frac{\mathbf{\dot{\gamma}}\times\mathbf{\ddot{\gamma}}}{\|\mathbf{\dot{\gamma}}\times\mathbf{\ddot{\gamma}}\|}$ and $N=B\times T$ we get equations \eqref{T}, \eqref{N}, \eqref{B}.
\hfill$\square$

The next statement gives us the relations between the focal curvatures of $\mathbf{\gamma}$ and the signed curvature of $\mathbf{\alpha}$, parameterized about an arbitrary parameter $t$.
\begin{theorem}\label{rel.euc.foc.cur.}
  Let $\mathbf{\alpha}=\mathbf{\alpha}(t),\, t\in I\subseteq\mathbb{R}$ be a regular $C^3$-plane curve in $\mathbb{E}^2$ and $K\neq0$ is the Euclidean signed curvature of $\mathbf{\alpha}$. Let $\mathbf{\gamma}(t)=\mathbf{\alpha}(t)+f(t).\vec{\mathbf{e}}_3,\,t\in I\subseteq\mathbb{R},\, f(t)\in C^2$ be the corresponding cylindrical curve over the right generalized cylinder with a base curve $\mathbf{\alpha}.$ If $c_1$ and $c_2$ are the focal curvatures of $\mathbf{\gamma}$, then
\begin{equation}\label{c1}
  c_1(t)=\displaystyle\frac{\sqrt{\dot{s}^{2}+\dot{f}^{2}}}
{\sqrt{\dot{s}^{6}K^{2}+\|\ddot{f}\dot{\alpha}-\dot{f}\ddot{\alpha}\|^{2}}},
\end{equation}
\begin{equation}\label{c2}
  c_2(t)=\displaystyle\frac{3(\dot{s}^{6}K^{2}+\|\ddot{f}\dot{\alpha}-\dot{f}\ddot{\alpha}\|^{2})
\frac{d}{dt}(\dot{s}^{2}+\dot{f}^{2})-(\dot{s}^{2}+\dot{f}^{2})
\frac{d}{dt}(\dot{s}^{6}K^{2}+\|\ddot{f}\dot{\alpha}-\dot{f}\ddot{\alpha}\|^{2})}
{2\sqrt{\dot{s}^{6}K^{2}+\|\ddot{f}\dot{\alpha}-\dot{f}\ddot{\alpha}\|^{2}}\,
(\dddot{f}\dot{s}^{3}K-<J(\ddot{f}\dot{\alpha}-\dot{f}\ddot{\alpha}),\dddot{\alpha}>)},
\end{equation}
where $\dot{s}$ is the arc-length prime of $\mathbf{\alpha}$ with respect to arbitrary parameter $t.$
\end{theorem}
{\it Proof:} From equation \eqref{GENC} and $c_1=\frac{1}{\varkappa}$ immediately follows equation \eqref{c1}.
After differentiating equation \eqref{c1} and replacing it along with equation \eqref{GENT} in  $c_2(t)=-\displaystyle\frac{\frac{d}{dt}\varkappa(t)}{\parallel \frac{d\mathbf{\gamma}(t)}{dt}\parallel \varkappa(t)^2 \tau(t)}=\displaystyle\frac{\frac{dc_1(t)}{dt}}{\parallel\frac{d\mathbf{\gamma}(t)}{dt}\parallel \tau(t)}$ we get \eqref{c2}.
\hfill$\square$

\begin{theorem}\label{focal via base curve}
  Let $\mathbf{\alpha}=\mathbf{\alpha}(t),\, t\in I\subseteq\mathbb{R}$ be a regular $C^3$-plane curve in $\mathbb{E}^2$ and $K\neq0$ is the Euclidean signed curvature of $\mathbf{\alpha}$. Let $\mathbf{\gamma}(t)=\mathbf{\alpha}(t)+f(t).\vec{\mathbf{e}}_3,\,t\in I\subseteq\mathbb{R},\, f(t)\in C^3$ be the corresponding cylindrical curve over the right generalized cylinder with base curve $\mathbf{\alpha}$ and $c_1, c_2$ are the focal curvatures of $\gamma$ defined by equations \eqref{c1} and \eqref{c2}. Then the focal curve of $\gamma$ has vector-parametric representation $C_{\gamma}(t)=\beta(t)+\tilde{f}(t).\vec{\mathbf{e}_{3}}$
\begin{equation}\label{beta}
  \beta(t)=\alpha(t)+\displaystyle\frac{c_{1}\left((\dot{s}^{2}+\dot{f}^{2})\ddot{\alpha}
-\frac{1}{2}\frac{d}{dt}(\dot{s}^{2}+\dot{f}^{2})\dot{\alpha}\right)-
c_{2}\sqrt{\dot{s}^{2}+\dot{f}^{2}}\, J(\ddot{f}\dot{\alpha}-\dot{f}\ddot{\alpha})}
{\sqrt{\dot{s}^{2}+\dot{f}^{2}}\,\sqrt{\dot{s}^{6}K^{2}+\|\ddot{f}\dot{\alpha}-\dot{f}\ddot{\alpha}\|^{2}}},
\end{equation}
\begin{equation}\label{f tilde}
 \tilde{f}(t)=f(t)+\displaystyle\frac{c_{1}\left(\ddot{f}\dot{s}^{2}-\dot{f}\frac{d}{dt}\left(\frac{\dot{s}^{2}}{2}\right)\right)+
c_{2}\sqrt{\dot{s}^{2}+\dot{f}^{2}}\,\dot{s}^{3}K}
{\sqrt{\dot{s}^{2}+\dot{f}^{2}}\sqrt{\dot{s}^{6}K^{2}+\|\ddot{f}\dot{\alpha}-\dot{f}\ddot{\alpha}\|^{2}}},
\end{equation}
where $\beta$ is \textbf{the generalized focal curve} of $\alpha$, $\dot{s}$ is the arc-length prime of $\mathbf{\alpha}$ with respect to arbitrary parameter $t.$
\end{theorem}
{\it Proof:} The proof follows immediately from the parametric representation of the focal curve given in equation~\eqref{C_gamma}$, i. e. \\\mathbf{ C_\gamma}(t)=\mathbf{\gamma}(t)+c_1(t) \mathbf{N}(t)+c_2(t) \mathbf{B}(t)$, using Theorem~\ref{frenet frame via base curve} and Theorem~\ref{rel.euc.foc.cur.}.
\hfill$\square$

\section{Applications}
\subsection{Closed toroidal curves}

Let us consider a torus with parametric equation
\begin{equation}\label{torus}
S_{2}(u,v)=((a+b\cos u)\cos v,(a+b\cos u)\sin v,b.\sin u),
\end{equation}
where $a,b=const, \, a>b>0, \, (u,v) \in D\subseteq \mathbb{E}^2$. The Cartesian equation of this surface of revolution obtained by the rotation of a circle in the $Oxz$ plane with center $a$ and radius $b$ around the $Oz$ axis is $(a-\sqrt{x^2+y^2})^2+z^2=b^2.$
Let $\gamma=S_{1}\cap S_{2}$ be the intersection curve of right generalized cylinder $S_1$ over a regular plane curve $\alpha$ and a torus $S_2$. Then the two equations that the scalar function $f: I \rightarrow\mathbb{R},\,I\subseteq\mathbb{R}$ satisfies are $f_{1,2}(t)=\pm\sqrt{b^{2}-(a-\sqrt{x^{2}(t)+y^{2}(t)})^{2}}$.

\begin{prop}\label{f_diferen}
  Since functions $x:I\rightarrow\mathbb{R}$ and $y:I\rightarrow\mathbb{R}$, $I\subseteq\mathbb{R}$ satisfy the inequalities $(a-b)^2 < x^2(t) + y^2(t) < (a+b)^2$ then functions $f_{1,2}:I\rightarrow\mathbb{R}$ given by $\pm\sqrt{b^{2}-(a-\sqrt{x^{2}(t)+y^{2}(t)})^{2}}$ have continuous partial derivatives of all orders.
\end{prop}
{\it Proof:}
The proof immediately follows from the condition $(a-\sqrt{x^{2}(t)+y^{2}(t)})^{2}<b^2$ which is true when functions $x:I\rightarrow\mathbb{R}$ and $y:I\rightarrow\mathbb{R}$ satisfied the inequalities $(a-b)^2 < x^2(t) + y^2(t) < (a+b)^2$.
\hfill$\square$

\begin{theorem}
 Let $\mathbf{\alpha}(t)=(x(t),y(t),0),\, t\in I\subseteq\mathbb{R}$ be a regular $C^3$-plane curve in $\mathbb{E}^2$ and $f: I\rightarrow\mathbb{R}$ is a differentiable scalar function $f(t)=\pm\sqrt{b^{2}-(a-\sqrt{x^{2}(t)+y^{2}(t)})^{2}}$ when $(a-b)^2 < x^2(t) + y^2(t) < (a+b)^2$. Then non-planar toroidal curve $\gamma : I\rightarrow \mathbb{E}^3$ given by $\mathbf{\gamma}(t)=\mathbf{\alpha}(t)+f(t).\vec{\mathbf{e}}_3$ is a curve of class $C^3$.
\end{theorem}
{\it Proof:}
The vector function $\mathbf{\gamma}(t)=\mathbf{\alpha}(t)+f(t).\vec{\mathbf{e}}_3=(x(t),y(t),f(t))$ is of class $C^3$ when its coordinate functions $x(t),y(t),f(t)$ are of class $C^3$. The fact that the functions $x(t),y(t)$ are continuously differentiable up to order $3$ immediately follows from the condition that $\alpha$ is a regular $C^3$-plane curve in $\mathbb{E}^2$. By Proposition \ref{f_diferen} the function $f(t)=\pm\sqrt{b^{2}-(a-\sqrt{x^{2}(t)+y^{2}(t)})^{2}}$ has continuous partial
derivatives of all orders therefore $\dot{f}(t)=\displaystyle\frac{df}{dt}(x(t),y(t))=\displaystyle\frac{df}{dx}\frac{dx}{dt}+\frac{df}{dy}\frac{dy}{dt}=f^{'}_{x}\dot{x}+f^{'}_{y}\dot{y}$, $\ddot{f}(t)=\displaystyle\frac{d^2f}{dt^2}(x(t),y(t)=f^{'}_{xx}\dot{x}^2+2f^{'}_{xy}\dot{x}\dot{y}+f^{'}_{yy}\dot{y}^{2}+f^{'}_{x}\ddot{x}+f^{'}_{y}\ddot{y}$ and  $\dddot{f}(t)=\displaystyle\frac{d^3f}{dt^3}(x(t),y(t))=f^{'}_{xxx}\dot{x}^3+(2f^{'}_{xyx}+f^{'}_{xxy})\dot{x}^2\dot{y}+(f^{'}_{yyx}+2f^{'}_{xyy})\dot{x}\dot{y}^{2}+f^{'}_{yyy}\dot{y}^3+
3(f^{'}_{xx}\dot{x}\ddot{x}+f^{'}_{xy}(\ddot{x}\dot{y}+\dot{x}\ddot{y})+f^{'}_{yy}\dot{y}\ddot{y})+f^{'}_{x}\dddot{x}+f^{'}_{y}\dddot{y}$  are continuously differentiable and then $f(t)\in C^3$.
\hfill$\square$

It can be shown that if the rulings of the generalized cylinder that passes though a cusp of the plane curve (generating the cylinder) are tangent to the meridians of the torus at the corresponding cusp, then the curve is connected.

\subsection{A helical curve over torus}

A curve that is orthogonal projection of a helix, wrapped into a torus (toroidal helix)
\[\gamma(t)=(\cos (t) (a+b \cos (n t)),\sin (t) (a+b \cos (n t)),b \sin (n t)),\, a,b,n=const, \, a,b,n>0,\]
 onto the Euclidean plane has parametric equation
$\alpha(t)=(\cos (t) (a+b \cos (n t)),\sin (t) (a+b \cos (n t)),0),\, t \in \mathbb{R}.$\\
 By Proposition \ref{focal via base curve} \textbf{the focal curve of toroidal helix} $C_{\gamma}(t)=\beta(t)+\tilde{f}(t).\vec{\mathbf{e}_{3}}$ is given by\\
$$x_{\beta}=\frac{4 a b n^2 \left(\cos (t) \left(a \left(n^2-1\right)-b \left(2 n^2+1\right) \cos (n t)\right)+3 b n \sin (t) \sin (n t)\right)}{\left(4 a^2 \left(n^2-1\right)-b^2 \left(8 n^4+13 n^2+3\right)\right) \cos (n t)+b \left(4 a \left(n^4-1\right)-4 a \left(2 n^2+1\right) \cos (2 n t)+b \left(n^2-1\right) \cos (3 n t)\right)},$$
$$y_{\beta}=\frac{4 a b n^2 \left(\sin (t) \left(a \left(n^2-1\right)-b \left(2 n^2+1\right) \cos (n t)\right)-3 b n \cos (t) \sin (n t)\right)}{\left(4 a^2 \left(n^2-1\right)-b^2 \left(8 n^4+13 n^2+3\right)\right) \cos (n t)+b \left(4 a \left(n^4-1\right)-4 a \left(2 n^2+1\right) \cos (2 n t)+b \left(n^2-1\right) \cos (3 n t)\right)},$$
$$\tilde{f}(t)=\frac{2 a \sin (n t) \left(-2 a^2 \left(n^2-1\right)+4 a b \left(2 n^2+1\right) \cos (n t)-b^2 \left(n^2-1\right) \cos (2 n t)+b^2 \left(11 n^2+1\right)\right)}{\left(4 a^2 \left(n^2-1\right)-b^2 \left(8 n^4+13 n^2+3\right)\right) \cos (n t)+b \left(4 a \left(n^4-1\right)-4 a \left(2 n^2+1\right) \cos (2 n t)+b \left(n^2-1\right) \cos (3 n t)\right)}.$$

\begin{figure}[!h]
\begin{center}
\begin{subfigure}{0.35\textwidth}
\includegraphics[width=\textwidth]{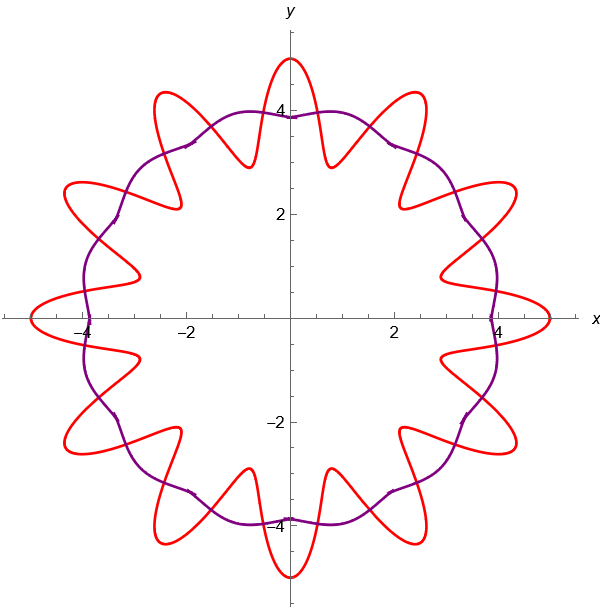}
\caption{\label{gen_foc_helix}curve $\alpha$ (in red), $\beta$ (in purple)}
 \end{subfigure}
%\hfill
 \begin{subfigure}{0.35\textwidth}
\includegraphics[width=\textwidth]{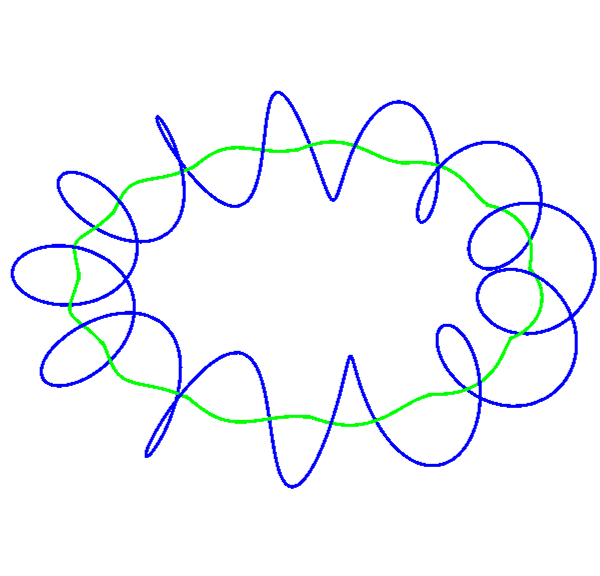}
\caption{\label{helix_fokal}curve $\gamma$ (in blue), $C_{\gamma}$ (in green)}
 \end{subfigure}
 \end{center}
\end{figure}
 The images of \textbf{orthogonal projection of toroidal helix} onto the Euclidean plane $\mathbf{\mathbf{\alpha}}$ (in red) and its corresponding \textbf{generalized focal curve} $\mathbf{\beta}$ (in purple) are in Figure~\ref{gen_foc_helix}. The images of \textbf{toroidal helix} $\mathbf{\gamma}$ (in blue) and its corresponding \textbf{focal curve} $\mathbf{C_\gamma}$ (in green) are in Figure~\ref{helix_fokal} for $a=4,b=1,n=12$.

\subsection{Epicycloids over torus}
Let us now consider the epicycloid given by
\[\alpha(t)=\left((r+R) \cos \left(\frac{r t}{R}\right)-r \cos \left(\frac{t (r+R)}{R}\right),(r+R) \sin \left(\frac{r t}{R}\right)-r \sin \left(\frac{t (r+R)}{R}\right),0\right),\, t\in [0,2k \pi], \, k=1,2,\ldots,\]
where $R$ is the radius of the larger circle centered at the origin around which the smaller circle with radius $r$ rolls.
We will take $R+2r\leq a+b$ and the larger circle to have radius $R=a-b,\, a>b$. Then the corresponding epicycloid $\alpha$ will be contained in the closed disk centered at the origin with radius $a+b$ (see \cite[p.4]{Cao2017}) and the rulings of the generalized cylinder that passes though a cusp of the epicycloid (generating the cylinder) will be tangent to the meridians of torus at the corresponding cusp and the toroidal epicycloid $\gamma$ will be connected.

The form of the curve of a epicycloid depends on the ratio $\frac{r}{R}=m$:

 - For $R=r\Leftrightarrow \frac{a}{2}\leq b < a$ the epicycloid $\alpha$ is called a \textbf{cardioid} and $\gamma$ is a \textbf{toroidal cardioid} with parametric equation
\[\gamma(t)=\left(r(2  \cos (t)- \cos (2 t)),r(2 \sin (t)- \sin (2 t)),\sqrt{(a-r)^2-(a-r\sqrt{5-4\cos (t)})^2}\right),\,t\in [0,2 \pi].\]
Let us take $R+2r<a+b \Leftrightarrow b>\frac{a}{2}$ and without losing of the generality let $b=\frac{3}{4}a \Leftrightarrow R=r=a-b=\frac{a}{4}=\frac{b}{3}$. Then $z_{\gamma}=r\sqrt{9-(4-\sqrt{5-4\cos (t)})^2}.$
%\[\gamma(t)=(r(2  \cos (t)- \cos (2 t)),r(2 \sin (t)- \sin (2 t)),r\sqrt{9-(4-\sqrt{5-4\cos (t)})^2}\,).\]
In this case $\gamma$ is non-planar when $t\neq {0,2\pi}$. The point $\gamma(0)=\gamma(2\pi)$ is a cusp of $\gamma$.

If $R+2r=a+b \Leftrightarrow a=2b \Leftrightarrow r=\frac{a}{2}=b$ then $z_{\gamma}=r\sqrt{1-(2-\sqrt{5-4\cos (t)})^2}$ and $\gamma$ is non-planar when $t\neq {0,\pi,2\pi}$ and the points $\gamma(0)=\gamma(2\pi)$, $\gamma(\pi)$ are cusps of $\gamma$.
\begin{figure}[!h]
\begin{center}
%\listoffigures
\begin{subfigure}{0.40\textwidth}
\includegraphics[width=\textwidth]{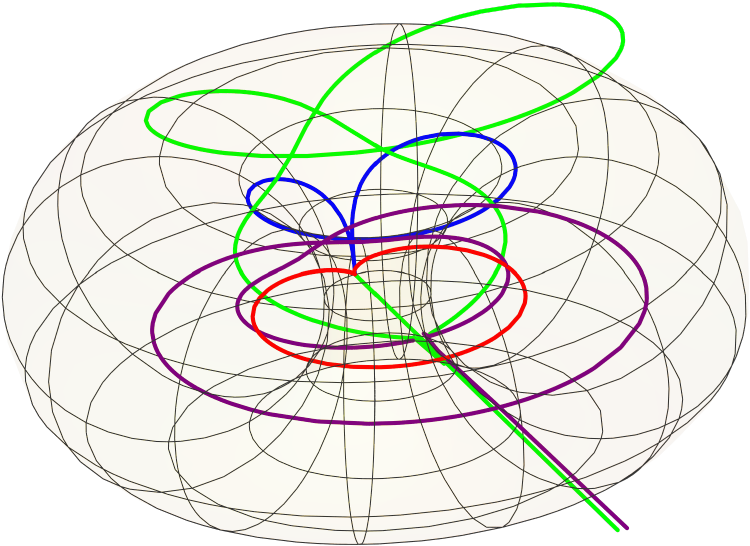}
\caption{\label{all_cardioid1}}
 \end{subfigure}
\begin{subfigure}{0.30\textwidth}
\includegraphics[width=\textwidth]{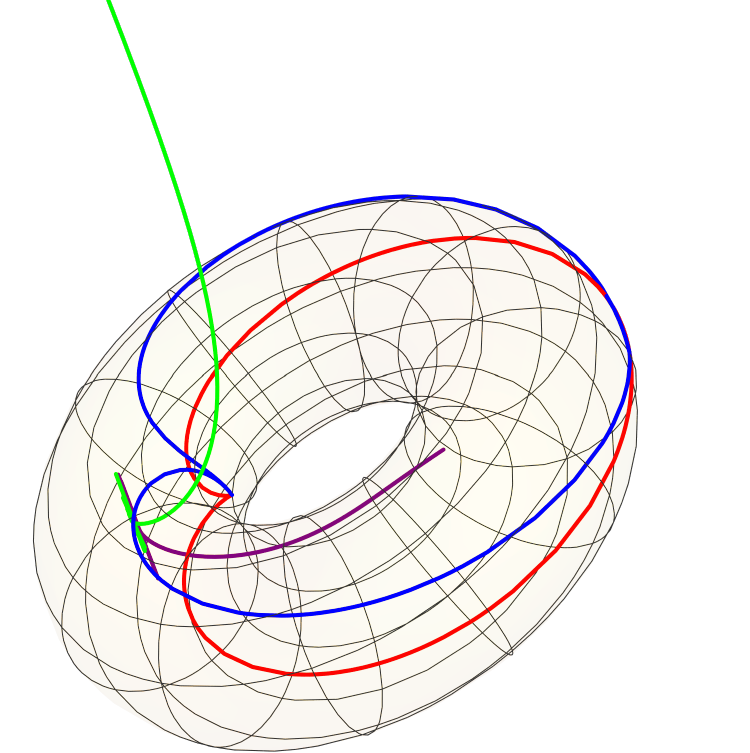}
\caption{\label{all_cardioid5}}
 \end{subfigure}
\end{center}
\end{figure}
The images of a \textbf{cardioid} $\mathbf{\alpha}$ (in red), corresponding \textbf{toroidal cardioid} $\mathbf{\gamma}$ (in blue) over torus, \textbf{the focal curve of toroidal cardioid} $\mathbf{C_\gamma}$ (in green) and \textbf{the generalized focal curve of the cardioid} $\mathbf{\beta}$ (in purple) are in Figure~\ref{all_cardioid1} (when $R+2r<a+b$) and Figure~\ref{all_cardioid5} (when $R+2r=a+b$) for $r=1$.

 - For $R=2r\Leftrightarrow \frac{a}{3}\leq b < a$ the epicycloid $\alpha$ is called a \textbf{nephroid} and $\gamma$ is a \textbf{toroidal nephroid} with parametric equation
\begin{center}
$\gamma(t)=\left(r\left(3  \cos \left(\frac{t}{2}\right)- \cos \left(\frac{3 t}{2}\right)\right),r\left(3 \sin \left(\frac{t}{2}\right)- \sin \left(\frac{3 t}{2}\right)\right),\sqrt{(a-2r)^2-(a-r\sqrt{10-6\cos (t)})^2}\right),\,t\in [0,4 \pi].$
\end{center}
Let us take $R+2r<a+b \Leftrightarrow b>\frac{a}{3}$ and without losing of the generality let $b=\frac{2}{3}a \Leftrightarrow R=2r=a-b=\frac{a}{3}=\frac{b}{2}$. Then $z_{\gamma}=r\sqrt{16-(6-\sqrt{10-6\cos (t)})^2}.$
%\[\gamma(t)=(r(2  \cos (t)- \cos (2 t)),r(2 \sin (t)- \sin (2 t)),r\sqrt{9-(4-\sqrt{5-4\cos (t)})^2}\,).\]
In this case $\gamma$ is non-planar when $t\neq {0,2\pi,4\pi}$. The points $\gamma(0)=\gamma(4\pi)$ and $\gamma(2\pi)$ are cusps of $\gamma$.
If $R+2r=a+b \Leftrightarrow a=3b=3r$ then $z_{\gamma}=r\sqrt{1-(3-\sqrt{10-6\cos (t)})^2}$ and $\gamma$ is non-planar when $t\neq {0,\pi,2\pi,3\pi,4\pi}$ and the points $\gamma(0)=\gamma(4\pi)$, $\gamma(\pi)$, $\gamma(2\pi)$, $\gamma(3\pi)$ are cusps of $\gamma$. In conclusion, we can say that if the epicycloid touches the circle of a torus with an equation $x^2+y^2=(a+b)^2$ then the cusps of the corresponding toroidal epicycloid are doubled.
\begin{figure}[!h]
\begin{center}
\begin{subfigure}{0.40\textwidth}
\includegraphics[width=\textwidth]{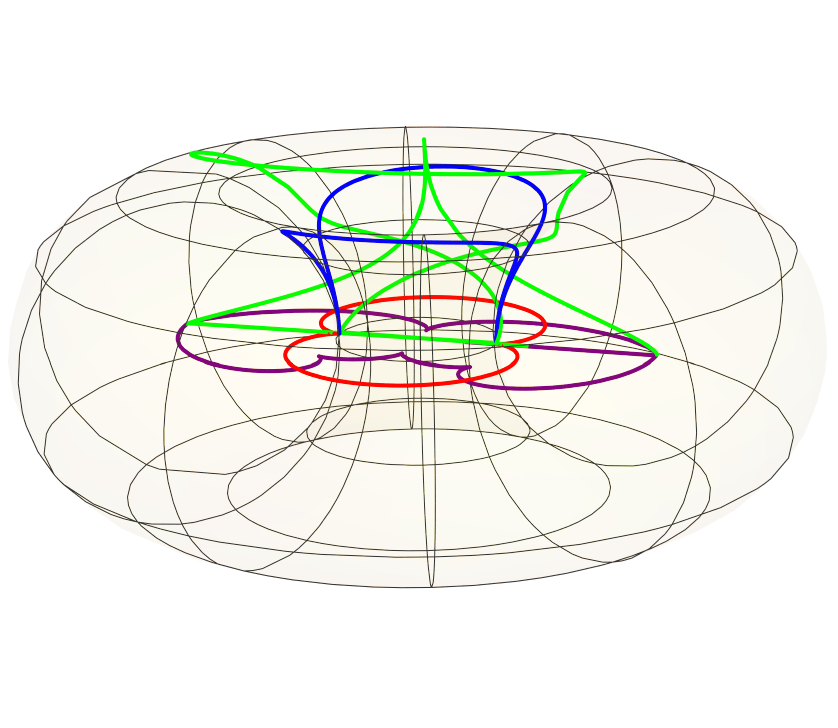}
\caption{\label{all_nephro1}}
\end{subfigure}
\begin{subfigure}{0.35\textwidth}
\includegraphics[width=\textwidth]{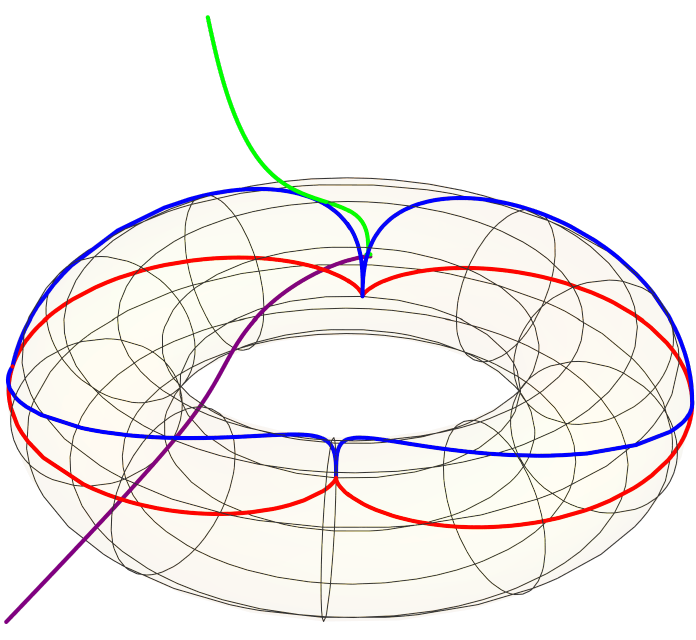}
\caption{\label{all_nephro2}}
\end{subfigure}
\end{center}
\end{figure}
The images of a \textbf{nephroid} $\mathbf{\alpha}$ (in red), corresponding \textbf{toroidal nephroid} $\mathbf{\gamma}$ (in blue) over torus, \textbf{the focal curve of toroidal nephroid} $\mathbf{C_\gamma}$ (in green) and \textbf{the generalized focal curve of the nephroid} $\mathbf{\beta}$ (in purple) are in Figure~\ref{all_nephro1} (when $R+2r<a+b$ for $r=1/2$) and Figure~\ref{all_nephro2} (when $R+2r=a+b$ for $r=1$).

\subsection{Hypocycloids over torus}
Let us now consider the hypocycloid given by
\[\alpha(t)=\left(r \cos \left(\frac{t (R-r)}{R}\right)+(R-r) \cos \left(\frac{r t}{R}\right),(R-r) \sin \left(\frac{r t}{R}\right)-r \sin \left(\frac{t (R-r)}{R}\right),0\right),\,t\in [0,2k\pi], k=1,2,..., \]
where $R$ is the radius of the larger circle centered at the origin on which the smaller circle with radius $r$ rolls internally.
We will take $|2r-R|\geq a-b,\, a>b$ and the larger circle to have radius $R\leq a+b$. Then the corresponding hypocycloid, will be contained in annulus $(a-b)^2\leq x^2+y^2\leq (a+b)^2$ (the region lying between two concentric circles centered at the origin with radii $a+b$ and $a-b$).
It is easy to see that from $|2r-R|\geq a-b,\, a>b$ and $R\leq a+b$ follows the inequality $R-r\leq2b$.
Then corresponding \textbf{toroidal hypocycloid} have representation $\gamma(t)=\alpha(t)\pm\sqrt{b^2-(a-\sqrt{(R-r)^2+2 r (R-r) \cos (t)+r^2}\,)^2}\,\,\vec{\mathbf{e}}_3.$ Let us take $R=a+b$.
\\The form of the curve of a hypocycloid depends on the ratio $\frac{r}{R}=m$:

- For $2R=3r$ the hypocycloid $\alpha$ is called a \textbf{deltoid} and $\gamma$ is a \textbf{toroidal deltoid} with parametric equation
\begin{center}
$\gamma(t)=\left(\frac{r}{2}\left(\cos (\frac{2t}{3})+2 \cos (\frac{t}{3})\right),\frac{r}{2}\left(\sin (\frac{2t}{3})- 2\sin (\frac{t}{3})\right),\sqrt{\left(\frac{3r}{2}-a\right)^2-\left(a-\frac{r}{2}\sqrt{5+4\cos (t)}\right)^2}\right), \, t\in [0,6\pi].$
\end{center}
If $|2r-R|>a-b,\, a>b$ then $\frac{a}{2}<b<a$ and if we take $b=\displaystyle\frac{2a}{3}$ then $R=\displaystyle\frac{3r}{2}=a+b=\displaystyle\frac{5a}{3}$.
In this case $z_\gamma=\frac{r}{10}\sqrt{36-(9-5\sqrt{5+4\cos (t)})^2}$ and
 $\gamma$ is non-planar when $t\neq {0,2\pi,4\pi,6\pi}$ and the points $\gamma(0)=\gamma(6\pi)$, $\gamma(2\pi)$, $\gamma(4\pi)$ are cusps of $\gamma$. If $2r-R=a-b,\, a>b$ then $b=\frac{a}{2}=\frac{r}{2}=\frac{R}{3}$. In this case $z_\gamma=\frac{r}{2}\sqrt{1-(2-\sqrt{5+4\cos (t)})^2}$
and $\gamma$ is non-planar when $t\neq {0,\pi,2\pi,3\pi,4\pi,5\pi,6\pi}$. The points $\gamma(0)=\gamma(6\pi)$, $\gamma(\pi)$, $\gamma(2\pi)$, $\gamma(3\pi)$, $\gamma(4\pi)$, $\gamma(5\pi)$ are cusps of $\gamma$.
\begin{figure}[!h]
\begin{center}
 \begin{subfigure}{0.45\textwidth}
\includegraphics[width=\textwidth]{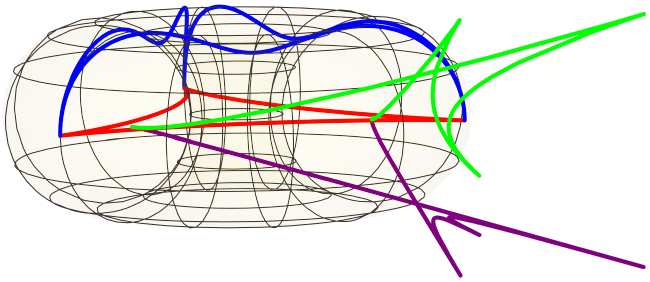}
\caption{\label{all_deltoid1}}
\end{subfigure}
% \hfill
\begin{subfigure}{0.45\textwidth}
\includegraphics[width=\textwidth]{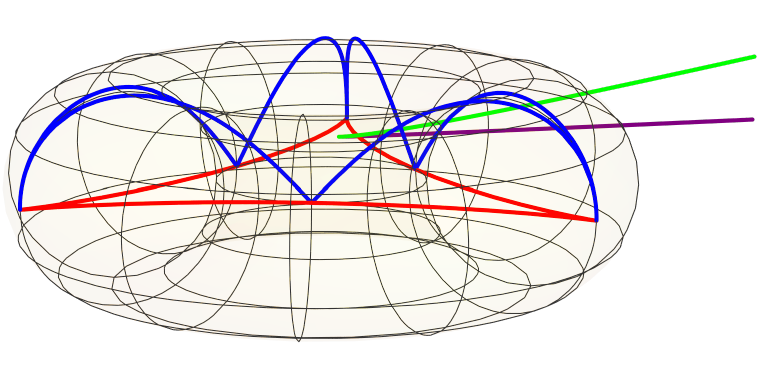}
\caption{\label{all_deltoid2}}
 \end{subfigure}
\end{center}
\end{figure}
 The images of \textbf{deltoid} $\mathbf{\alpha}$ (in red), the corresponding \textbf{toroidal deltoid} $\mathbf{\gamma}$ (in blue) over torus, \textbf{the focal curve of toroidal deltoid} $\mathbf{C_\gamma}$ (in green) and \textbf{the generalized focal curve of the deltoid} $\mathbf{\beta}$ (in purple) are in Figure~\ref{all_deltoid1} (when $|2r-R|>a-b,\, a>b$ for $r=10$) and Figure~\ref{all_deltoid2} (when $|2r-R|=a-b,\, a>b$ for $r=2$).

- For $R=4r$ the hypocycloid $\alpha$ is called an \textbf{astroid} and $\gamma$ is a \textbf{\textbf{toroidal astroid}} with parametric equation
\begin{center}
$\gamma(t)=\left(r(3  \cos (\frac{t}{4})+ \cos (\frac{3 t}{4})),r(3 \sin (\frac{t}{4})- \sin (\frac{3 t}{4})),\sqrt{(4r-a)^2-(a-r\sqrt{10+6\cos (t)})^2}\right), t\in[0,8\pi].$
\end{center}
If $|2r-R|>a-b,\, a>b$ then $\frac{a}{3}<b<a$ and if we take $b=\displaystyle\frac{a}{2}$ then $R=4r=a+b=\displaystyle\frac{3a}{2}$.
In this case $$\gamma(t)=\left(4r\cos^3\left(\frac{t}{4}\right),4r\sin^3\left(\frac{t}{4}\right),\frac{r}{3}\sqrt{16-(8-3\sqrt{10+6\cos (t)})^2}\right)$$ and $\gamma$ is non-planar when $t\neq {0,2\pi,4\pi,6\pi,8\pi}$ and the points $\gamma(0)=\gamma(8\pi)$, $\gamma(2\pi)$, $\gamma(4\pi)$, $\gamma(6\pi)$ are cusps of $\gamma$. If $|2r-R|=a-b,\, a>b$ then $b=r=\frac{a}{3}=\frac{R}{4}$.
In this case $z_\gamma=r\sqrt{1-(3-\sqrt{10+6\cos (t)})^2}$ and $\gamma$ is non-planar when $t\neq {0,\pi,2\pi,3\pi,4\pi,5\pi,6\pi,7\pi,8\pi}$ and the points $\gamma(0)=\gamma(8\pi)$, $\gamma(\pi)$, $\gamma(2\pi)$, $\gamma(3\pi)$, $\gamma(4\pi)$, $\gamma(5\pi)$, $\gamma(6\pi)$, $\gamma(7\pi)$ are cusps of $\gamma$.
In conclusion, we can say that if the hypocycloid touches the circle of a torus with an equation $x^2+y^2=(a-b)^2$ then the cusps of the corresponding toroidal hypocycloid are doubled.
\begin{figure}[!h]
\begin{center}
\begin{subfigure}{0.45\textwidth}
\includegraphics[width=\textwidth]{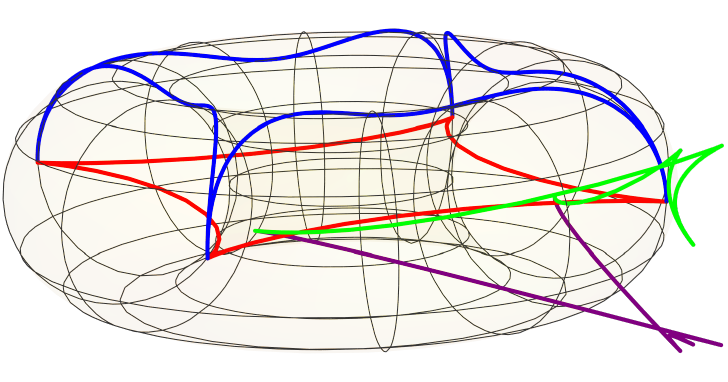}
\caption{\label{all_astroid1}}
 \end{subfigure}
% \hfill
  \begin{subfigure}{0.50\textwidth}
\includegraphics[width=\textwidth]{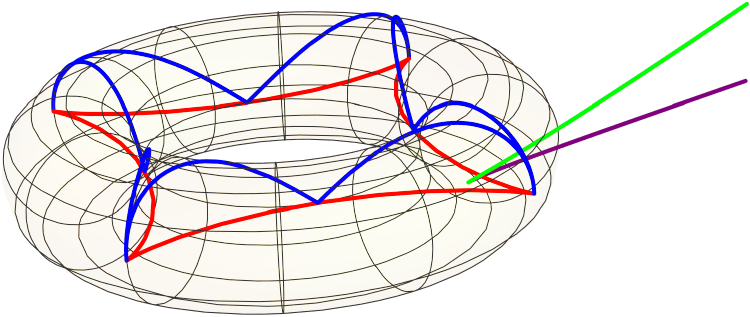}
\caption{\label{all_astroid2}}
 \end{subfigure}
\end{center}
\end{figure}
The images of a \textbf{astroid} $\mathbf{\alpha}$ (in red), the corresponding \textbf{toroidal astroid} $\mathbf{\gamma}$ (in blue) over torus, \textbf{the focal curve of toroidal astroid} $\mathbf{C_\gamma}$ (in green) and \textbf{the generalized focal curve of the astroid} $\mathbf{\beta}$ (in purple) are in Figure~\ref{all_astroid1} (when $|2r-R|>a-b,\, a>b$ for $r=3/4$) and Figure~\ref{all_astroid2} (when $|2r-R|=a-b,\, a>b$ for $r=1$).

\section{Cross-referencing}

\section{Conclusion}
In this paper, we discuss relations between differential geometric invariants of non-planar space curve on the right generalized cylinder over regular base curve and differential geometric invariants of corresponding base curve. Based on these relations, a method for obtaining a new space curve from a given plane curve parameterized about an arbitrary parameter is presented. At first, we define a non-planar space
curve on the right generalized cylinder whose base curve is the considered plane curve, parameterized about an arbitrary parameter.
Later on, we examine the focal curve of the obtained cylindrical curve which is also a non-planar curve. The Frenet-Seret system
of the cylindrical curve and its focal curve are expressed in terms of the signed curvature of the abovementioned plane curve and
its derivatives. Finally, we obtained the parametric representation of orthogonal projection of the focal curve onto the Euclidean
plane via aforementioned plane curve and its derivatives. That curve is called generalized focal curve of a plane curve. The
proposed method is demonstrated for several closed plane curves used in engineering practice. These curves include: epicycloid,
hypocycloid and a curve that is orthogonal projection of toroidal helix onto the Euclidean plane. Moreover, depending on the type of closed plane curves, conditions are derived where their corresponding space curves on the torus are also closed.

\section*{Author contributions}
All authors contributed equally to the preparation of the study: Conceptualization, Methodology, Investigation, Writing-original draft, Writing-review \& editing.
%, and are listed in alphabetical order.
All authors have read and approved the final version of the manuscript for publication.

\section*{Use of AI tools declaration}
The authors declare they have not used Artificial Intelligence (AI) tools in the creation of this article.

%\section*{Acknowledgments (All sources of funding of the study must be disclosed)}
\begin{acknowledgments}
The first author is partially supported by Scientific Research Grant RD-- 08--104/30.01.2024 of Konstantin Preslavsky University of Shumen.
%The second author is partially financed by Scientific Research Grant RD--08--/30.01.2024 of Konstantin Preslavsky University of Shumen.
The third author is partially supported by Scientific Research Grant RD-- 08--103/30.01.2024 of Konstantin Preslavsky University of Shumen.
\end{acknowledgments}

\section*{Conflict of interest}
The authors declare no conflict of interest.

\nocite{*}
\begin{footnotesize}
%\bibliography{aipsamp}% Produces the bibliography via BibTeX.

\begin{thebibliography}{7}

\bibitem{article1}
Encheva, R.P., Georgiev, G.H., Curves on the shape spere, Result. Math. {\bf 44} (2003), 279--288.

\bibitem{article2}
Encheva, R.P., Georgiev, G.H., Shapes of space curves, J. Geom. Graph. {\bf 7} (2003), 145--155.

\bibitem{article3}
Encheva, R.P., Georgiev, G.H., Similar Frenet Curves, Result. Math., {\bf 55} (3-4) (2009), 359--372.
http://www.springerlink.com/content/x6145691n7437r17/

\bibitem{Georgiev2015}
Georgiev, G.H., Encheva, R.P., Dinkova, Cv.L., Geometry of cylindrical curves over plane curves, Applied Mathematical Sciences,{\bf 113} (9) (2015), 5637--5649. http://dx.doi.org/10.12988/ams.2015.56456

\bibitem{article5}
Arroyo, J., Garay, O.J., Menca, J.J., When is a peroidic function the curvature of a closed plane curve?, Am. Math. Mon. {\bf 115} (5), (2008), 405--414.

\bibitem{article6}
Izumiya, S., Takeuchi, N., Special curves and ruled surfaces, Cotributions to Algebra and Geometry, {\bf 44} (2003), 203--212.

\bibitem{Gray2006}
Gray, A., Abbena, E., Salamon, S., Modern Differential Geometry of Curves and Surfaces, Chapman Hall/CRC, (2006).

\bibitem{Do Carmo2016}
Do Carmo, M.P., Differential Geometry of Curves and Surfaces: Revised and Updated Second Edition, Dover Publications, INC, (2016)

\bibitem{Cao2017}
Cao, C., Fletcher, A., Ye, Z. (2017). Epicycloids and Blaschke products. Journal of Difference Equations and Applications, 23(9), 1584–1596. https://doi.org/10.1080/10236198.2017.1341499

\end{thebibliography}

\end{footnotesize}
\end{document}